\documentclass[12pt]{article}

\usepackage{amsmath}
\usepackage{amscd,amsthm}
\usepackage{psfig}
\usepackage{graphicx}
\usepackage{epsfig, graphicx,subfigure}
\font\bb=msbm10 at 12pt
\def\rR{\hbox{\bb R}}

\newtheorem{lemma}{Lemma}

\newcommand{\beq}{\begin{eqnarray}}
\newcommand{\eeq}{\end{eqnarray}}
\newcommand{\beqq}{\begin{eqnarray*}}
\newcommand{\eeqq}{\end{eqnarray*}}

\theoremstyle{plain}
\newtheorem{Thm}{Theorem}[section]
\newtheorem{Prop}[Thm]{Proposition}
\newtheorem{Cor}[Thm]{Corollary}
\newtheorem{Lemm}[Thm]{Lemma}
\newtheorem{AThm}{Boundary Theorem}

\newtheorem{Corr}{Corollary}

\newtheorem{BThm}{Linear
Kick Theorem}

\newtheorem{LogThm}{Logarithmic
Kick Theorem}

\newtheorem{HookThm}{Hook Theorem}

\theoremstyle{remark}
\newtheorem*{Defn%
}{\bf Definition}
\newtheorem*{Term%
}{\bf Terminology}
\newtheorem*{Rmk%
}{\bf Remark}
\newtheorem*{Notn%
}{\bf Notation}
\newcommand{\<}{\langle}
\renewcommand{\>}{\rangle}
\newcommand{\abs}[1]{\lvert#1\rvert}

\newcommand{\hardhyphen }%
{\nobreakdash-\hspace{0pt}}
\newcommand{\dlim}{\displaystyle\lim}

\newtheorem{theorem}{Theorem}

\newtheorem{conj}{conjecture}
\newtheorem{prop}{Proposition}

\font\bb=msbm10 at 12pt

\def\rR{\hbox{\bb R}}

\def\QED{\quad\hbox{\hskip 4pt\vrule width 5pt height 6pt depth 1.5pt}}

\begin{document}

\title{The Boundary between Compact and Noncompact Complete Riemann Manifolds}
\author{D. Holcman\thanks{Weizmann Institute of
Science,  Department of Mathematics, Rehovot 76100, Israel.email:david.holcman@weizmann.ac.il} \and C. Pugh \thanks{Mathematics Department, University of California, Berkeley California, 94720, U.S.A. email:pugh@math.berkeley.edu}}
\date{}
\maketitle
\begin{abstract} In 1941 Sumner
Myers proved that if the Ricci
curvature of a complete Riemann
manifold has a positive infimum
then the manifold is compact and
its diameter is bounded in terms
of the infimum. Subsequently the
curvature hypothesis has been
weakened, and in this paper we
weaken it further in an attempt to
find the ultimate, sharp result.
\end{abstract}



\section{Introduction}

Myers' Theorem  \cite{Myers}
states that a complete  Riemann
manifold
$(M,g)$ of dimension
$n \geq 2$ is compact if its
Ricci curvature is uniformly
positive, and furthermore it
has  diameter
$
\leq
\pi /\sqrt{C}$ if its Ricci
curvature satisfies
\begin{eqnarray}
\operatorname{Ric}_p   \geq
(n-1)C
\end{eqnarray} everywhere on
$M$,
$C$ being a positive constant.
(Here and below we adopt the
shorthand that
$\operatorname{Ric}_p \geq
c$ means that for all  $X \in
T_pM$,
$$
\operatorname{Ric}_p(g, X, X)
\geq c \<X, X\>_p ,
$$ where $\<\;\; , \;\; \>_p$ is the
$g$-inner product on $T_pM$).
Later,   asymptotic  conditions
on the curvature were found that
still imply compactness,
\cite{CGT}, \cite{Kupka},
\cite{BF}, although definitive
conditions remain unknown.
The idea is to fix an origin $O
\in  M$ and study the curvature
along geodesics   emanating
from $O$.  One always assumes
  the curvature is positive, but
permits it to decay to zero far
from
$O$.

To be more specific, we set
$$
\operatorname{Ric}(r) = \inf
\{\operatorname{Ric}_p : p =
\exp _O(v) \textrm{ and }
\abs{v} = r \},
$$
and assume throughout that
$\operatorname{Ric}(r) > 0$.
Hypotheses that imply
compactness and give a
diameter estimate are:
\begin{itemize}

   \item[(a)]
(Cheeger-Gromov-Taylor
\cite{CGT})
  For some $\nu  >
0$, some $r_0 > 0$, and all $r
\geq r_0$,
$$
\operatorname{Ric}(r) \geq
\dfrac{n-1}{4}\Big(\frac{1
+4\nu ^2}{r^2} \Big).
$$

   \item[(b)] (Cheeger-Gromov-Taylor
\cite{CGT})
For some $\nu  >
0$, some $r_0 > 1$, and all $r
\geq r_0$,
$$
\operatorname{Ric}(r) \geq
\dfrac{n-1}{4}\Big(
\frac{1}{r^2} + \frac{1
+4\nu ^2}{(r \ln r)^2} \Big).
$$

   \item[(c)]  (Boju-Funar
\cite{BF}) For some $\nu  >
0$, some integer $k $,
some
$r_0 > e_k$, and all
$r
\geq r_0$,
$$
\operatorname{Ric}(r) \geq
\dfrac{n-1}{4}\Big(
\frac{1}{r^2} +
\frac{1}{(r \ln r)^2} +  \dots +
\frac{1
+4\nu ^2}{(r \ln ( r) \ln \ln (r)
\cdots \ln ^k(r))^2}
\Big),
$$
where   $\ln ^k$
is the
$k^{\textrm{th}}$ iterated
  logarithm,
$\ln
^k(r) =
\ln \circ
\ln \circ   \dots \circ \ln(r)$, and
$\ln ^k(e_k) = 0$.

\end{itemize}
It is natural to   set
$\ln ^0 (r) = r$.  Then  (a)
and (b) are (c) with
$k = 0$ and $k = 1$.  The
diameter estimates on $M$
involve
$\nu$, $r_0$, and $k$.  When
$k=0$, one has
$\operatorname{diam}(M) \leq
e^{\pi /\nu}r_0$.  See
\cite{CGT} and \cite{BF}.

\begin{Rmk}
In \cite{Kupka}, Dekster and Kupka
prove that the estimate in (a) is
sharp.
\end{Rmk}

\begin{Rmk}
In terms of decay rates, these
results are nearly optimal.  For
example, there exist
noncompact complete
manifolds whose Ricci
curvature satisfies a
Boju-Funar equality with $\nu
= 0$,
\begin{equation}
\begin{split}
\label{e:BFE}
  \operatorname{Ric}(r) =
\dfrac{n-1}{4}\Big(
\frac{1}{r^2} +
   \dots +
\frac{1
}{(r \ln ( r) \ln \ln (r)
\cdots \ln ^k(r))^2}
\Big),
\end{split}
\end{equation}
   See \cite{BF}.
\end{Rmk}

\begin{Rmk}
Dekster and Kupka consider
also the \emph{sectional}
curvature $K$, and show that
all noncompact complete
Riemann manifolds of positive
curvature satisfy
$$
\liminf_{r \rightarrow \infty}
k(r) r^2 \leq \frac{1}{4},
$$
where $k(r) = \inf \{ K_p :  p =
\exp _O(v) \textrm{ and }
\abs{v} = r \}$.  The constant
$1/4$ is sharp, \cite{Kupka}.
\end{Rmk}

Our first   results, the Kick
Theorems,  state that
asymptotic estimates are not the
only way to guarantee
compactness.  Instead of
curvature that decays to zero at
a   positive rate ($\nu >
0$) as the radius tends to infinity,
it is enough that in addition to
Ricci curvature obeying
(\ref{e:BFE}), there is a certain
amount of extra curvature  on a
finite  shell $\{ p =
\exp _O(v) : a \leq
\abs{v} \leq  b \}$.  We refer to
the extra curvature as  a ``kick.''
See Sections ~\ref{s:kick} and
~\ref{s:log} for details.

\begin{Rmk} None of these
conditions is truly optimal;  in
Section~\ref{s:embedded} we
show that a surface
approximating the
    capped cylinder has the
property that every more curved
surface is compact, but this is
implied by none of the
asymptotic or kick conditions.
\end{Rmk}

Nevertheless, it is tempting to
postulate some kind of a
boundary in the space of
positive Ricci curvature functions
with all compact manifolds on one
side and all non-compact complete
manifolds on the other.\footnote{
\linespread{1.0}
\footnotesize Readers familiar
with Walter Rudin's text,
\emph{Principles of
Mathematical Analysis}, will
recognize this phrase, in which
Rudin asserts that there is no such
boundary dividing convergent
and divergent series.  A
difference between series and
curvature functions is that local
perturbations have no global
effect on series, while for
curvature functions this is not so.
Perturbation of a finite number
of terms in a series does not
change convergence, but a
compactly supported
perturbation of curvature can
affect the manifold's topology at
a long distance from the
perturbation's support, and
hence such a curvature boundary
is not unreasonable.} A
manifestation of such a
boundary would be a topology on
the space
$\mathcal{R}$ of Ricci
curvature functions which are
defined on a fixed tangent space
$T_OM$, and a  closed subset
$\mathcal{R}_0 \subset
\mathcal{R}$
    such that    through each
$R_0 \in \mathcal{R}_0$   there
is a
     curve
$R_t$ of Ricci curvature
functions, and
\begin{itemize}

    \item[(a)] If $t > 0$ and $M$
has   Ricci curvature
$ R_t$  then it  is   compact.

    \item[(b)] If $t \leq 0$ and
   $M$ has Ricci curvature $ R_t$
then it  is non-compact.

    \item[(c)]

This transverse, single-point
crossing from non-compact to
compact persists for all nearby
curves
$\widetilde{R }_t$.
\end{itemize} In
Section~\ref{s:planar} we
establish this kind of result for
planar curves, finding a
boundary between
embedded curves and
nonembedded immersed curves;
in Section~\ref{s:embedded} we
pass to  surfaces embedded in
$3$-space.

Our second result partially
identifies   the boundary
$\mathcal{R}_0$
    postulated   above.   We call a
function $b : [0, \infty)
\rightarrow  (0, \infty)$ an
\textbf{SL-bifurcator} if  the
solution to the Sturm-Liouville
equation
$$
w^{\prime}{}^{\prime} + b(r)w
= 0, \quad \quad w(0) = 0, \quad
w^{\prime}(0) = 1
$$
is monotone and bounded.
   An example is
$$
b(r) = \frac{2r}{(1 + r^2)^2
\arctan r}.
$$
See
Sections~\ref{s:boundary}
and \ref{s:embedded} for more
on SL-bifurcators.

  \begin{Defn}
A function $f(x)$
\textbf{exceeds} a function
$g(x)$, if for all $x$, $f(x) \geq
g(x)$, and for some $x$, $f(x) >
g(x)$.
\end{Defn}

\begin{AThm}
\label{t: } Let $M$ be a complete
Riemann manifold with positive
Ricci curvature, and let
$b(r)$ be an SL-bifurcator.
	\begin{itemize}

    \item[(a)] $M$ is noncompact
if   for all
$r
\geq 0$,
$$
\sup \{ \operatorname{Ric}_p :
p = \exp (v) \textrm{ and }
\abs{v} = r  \}  \leq
b(r).
$$

    \item[(b)]  $M$ is compact if
$\operatorname{Ric}(r)$
exceeds $b(r)$.   (As above,
$\operatorname{Ric}(r)$
denotes the infimum of $
\operatorname{Ric}_p $ such
that $p =
\exp _O(v) \textrm{ and }
\abs{v} = r $.)
\end{itemize}

\end{AThm}

\begin{Corr}
\label{c:boundary}
SL-bifurcators distinguish
compact and noncompact
complete Riemann manifolds
with positive Ricci curvature.
\end{Corr}
See
Section~\ref{s:boundary} for the
  simple proofs.

\section{Linear Kick
   }
\label{s:kick}

In this section we deal with the
kick condition when
$k = 0$.  Thus we assume that for
all $r
\geq r_0 > 0$,
\begin{equation}
\begin{split}
\label{e:base}
  \operatorname{Ric}(r) \geq
\frac{n-1}{4 }\Big(
\frac{1}{r^2}\Big),
\end{split}
\end{equation}
and we find a sufficient amount
of extra curvature on a shell $\{
p = \exp _O(v) : a \leq \abs{v}
\leq  b\}$ (with $r_0 \leq  a$) that
implies
$M$ is compact. Define
$\lambda
$ to be the smallest positive root
of the equation
\begin{equation}
\begin{split}
\label{e:nu}
    \cot (\lambda  \ln (b/a)) =
    \lambda  \ln (a/r_0).
\end{split}
\end{equation}

Note that $\lambda  = \lambda
(a,b, r_0)$ exists, lies in $(0, \pi /2
\ln (b/a)]$,  and is unique.  For,
as
$\mu \ln (b/a)$ varies from
$0$ to $\pi /2 $,  its cotangent
decreases monotonically from
$\infty$ to $0$, while  $\mu \ln
(a/ r_0)$ is non-decreasing.

\begin{BThm}
\label{t:our} In addition to
(\ref{e:base}), assume that for all
$r
\in [a, b]$, we have
$$
\operatorname{Ric}(r) >
\frac{n-1}{4}\Big(
\frac{1 + 4\lambda  ^2}{r^2}
  \Big),
$$
where $\lambda = \lambda
(a, b, r_0)$.  Then
$M$ is compact.
\end{BThm}


\begin{Rmk} As an example of
the kick we can take
$a = e$ and
$b = e^2$.  Then $\lambda $ is
approximately $.46$.  Also, if the
interval
$[a,b]$ is small, in the sense that
$b-a=\epsilon$, it  is not hard to
check that   a kick sufficient  for
compactness  increases like
$\epsilon^{-1/2}$ as $\epsilon
\rightarrow 0$.
\end{Rmk}

The proof   is
based on analyzing a kicked
Sturm-Liouville equation
\begin{eqnarray}
\label{e:KSL}
y^{\prime}{}^{\prime} +
\frac{1}{4}
\Big(\frac{1 + 4 \mu^2
\chi_{[a, b]} (r)}{r^2}
\Big) y =   0.
\end{eqnarray}

\begin{Lemm}
\label{l:SL} If
$\lambda  = \lambda  (a, b,
r_0)$ is determined as in
(\ref{e:nu}) and  if
$\mu  > \lambda  $ then the
solution
$y(r)$ of (\ref{e:KSL}) with
initial conditions
$y(r_0) = 0$,
$y^{\prime}(r_0) > 0$,
necessarily vanishes  at some
$r >
   r_0$.
\end{Lemm}

\begin{proof} For any constants
$c, k >  0$, the   function
$w(r) = cy(r/k)$ satisfies
(\ref{e:KSL}) and has initial
conditions
$$ w(kr_0) = 0 \qquad
w^{\prime}(kr_0) =
\frac{cy^{\prime}(kr_0)}{k}
    > 0.
$$ Taking
$k = 1/r_0$ and
$c = k/y^{\prime}(r_0)$, we can
    assume without loss of
generality that
$r_0 = 1$ and $y^{\prime}(1) =
1$. We do so.

For $r  \in [a, b]$, the solution of
(\ref{e:KSL})   is of the form
\begin{eqnarray}
    \label{e:se1} y(r) =A r^{1/2}
\cos(\mu \ln r) +B r^{1/2}
\sin(\mu
\ln r)
\end{eqnarray} where $A,B$ are
constants. For $\mu =0$, the
solution degenerates as
\begin{eqnarray}
\label{se2} y(r) =
Ar^{1/2}+Br^{1/2}
\ln r ,
\end{eqnarray} which can be seen
by replacing $B$ with  $B/ \mu $
in (\ref{e:se1}) and letting $\mu $
tend to zero. Matching initial
conditions at
$r = 1$, $r = a$, and $r= b$ gives
$$ y(r) \, = \,
\begin{cases}
    r^{1/2} \ln(r) &\textrm{ if } 1
\leq  r
\leq a \\
    r^{1/2}\Big( \ln(a)  \cos(\mu
\ln(r/a) ) +
\displaystyle
\frac{1}{\mu}\sin(\mu
\ln(r/a) \Big) &\textrm{ if }  a
\leq r
\leq b \\
    (r/b)^{1/2}( \alpha +
\beta  \ln (r/b) ) &\textrm{ if } b
\leq r <
\infty
\end{cases}
$$ where $\alpha , \beta $ are
constants
\begin{eqnarray}
\alpha  = b^{1/2} \Big( \ln(a)
\cos (\mu  \ln (b/a)) +
\frac{1}{\mu}
\sin (\mu  \ln (b/a) ) \Big)
\\
\beta  =b^{1/2}
\Big(
    \cos(\mu \ln(b/a))
    - \mu \ln(a)
\sin(\mu \ln (b/a))\Big)
\end{eqnarray}

    By the Sturm Comparison
Theorem, a second zero of
$y(r)$, if it exists, is a monotone
decreasing function of $\mu $.
Thus it is no loss of generality to
assume that
$\mu  -\lambda $ is small.  Since
$\lambda $ is the smallest
positive root of
$$
\cot (\lambda  \ln (b/a) ) =
\lambda
\ln a
\geq  0,
$$ and since $\lambda  \ln (b/a)
\leq
\pi /2$, we can assume $\mu  -
\lambda $ so small that
  $\mu
\ln (b/a) < \pi
$ and
$$
\cot (\mu \ln (b/a)) >
\frac{-1}{\mu \ln a}.
$$
Since the cotangent is
monotone decreasing on $(0, \pi
)$, this   implies that
$$
\cot (\mu  \ln (r/a)) >
\frac{-1}{\mu \ln a}
$$
for $a \leq  r \leq  b$, and
hence that
$y(r) > 0$ on
$[a, b]$.  For the same reasons,
$$
\cot (\mu  \ln (b/a)) < \cot
(\lambda
\ln (b/a)) = \lambda  \ln a  \leq
\mu
\ln a,
$$
which implies that $\beta  <
0$. But $ y(b) > 0$ and $\beta  <
0$ implies that
$y(r) = 0$ for some $r > b$.
\end{proof}

\begin{proof}[\bf Proof of the
Linear Kick  Theorem]    We
must show that $M$ is compact.
By the Hopf-Rinow Theorem, it
suffices to show that every
geodesic
   through
$O$ contains a pair of conjugate
points.  For if $M$ is not
compact then it contains an
everywhere
   distance minimizing geodesic
$\gamma
$ from
$O$ to infinity ($M$ is
complete), and this is contrary
to   conjugate pairs on
$\gamma $.  See
\cite{Kobayashi}.

By the assumption on the Ricci
curvature and compactness of
$[a, b]$, there exists $\mu $ such
that
$$
\operatorname{Ric}\, (\gamma
^{\prime}(r), \gamma
^{\prime}(r)) >
\mu  >
\lambda  =
\lambda  (a, b, r_0)
$$
for $r  \in [a, b]$.  Fix such a
$\mu
$, and let
$y(r)$ be a solution
    of the kicked Sturm-Liouville
equation (\ref{e:KSL}) with
initial conditions $y(r_0) = 0$,
$y^{\prime}(r_0) > 0$.  By
Lemma~\ref{l:SL}, $y(r_1) = 0$
for some $  r_1 > r_0$.

Following Myers' use of the Index
Theorem,  this gives a pair of
conjugate points on $\gamma
$.   We recapitulate his proof.

Choose an orthonormal basis
$\{e_1,
\dots  , e_n\}$ of
$T_OM$ with $e_n = \gamma
^{\prime}(0)$, and let $E_1 (r),
\dots , E_{n-1}(r)$ be the
corresponding parallel vector
fields along $\gamma $.  Define
vector fields along $\gamma $,
$$ X_j(r)  =  y(r) E_j(r).
$$ We will check that
\begin{equation}
\begin{split}
\label{e:sum}
    \sum_{j=1}^{n-1} I(X_j, X_j)
< 0
\end{split}
\end{equation} where
$$
I(X,X) = \int_{r_0}^{r_1} \Big(
\abs{\nabla_t X}^2 -
\< R(X, \gamma ^{\prime}), X
\> \Big)\, dr
$$
  is the index of a vector field
$X$ along $\gamma $.  ($R$ is
the curvature tensor.) To verify
(\ref{e:sum}) we evaluate the
Ricci curvature hypothesis on
$(\gamma ^{\prime}(r), \gamma
^{\prime}(r))$  as
\begin{equation*}
\begin{split} &
\sum \<R(E_j(r), \gamma
^{\prime}(r))\gamma
^{\prime}(r), E_j (r) \>
\\ &=
\operatorname{Ric}\, (\gamma
^{\prime}(r), \gamma
^{\prime}(r)) >\frac{n-1}{4}
\Big( \frac{1
+ 4 \mu ^2\chi_{[a,
b]}(r)
}{r^2} \Big).
\end{split}
\end{equation*} Then we write
\begin{equation*}
\begin{split}
\sum I_j  &= \sum \int
\<X_j^{\prime},
X_j^{\prime}\> -
\<R(X_j,
\gamma ^{\prime}) \gamma
^{\prime}, X_j\> \, dr
\\ &= \sum \int (y^{\prime})^2 -
\<R(E_j, \gamma
^{\prime})\gamma ^{\prime},
E_j\>y^2 \, dr
\\ & = (n-1) \int (y^{\prime})^2
\, dr -
\int
\sum
\<R(E_j, \gamma
^{\prime})\gamma ^{\prime},
E_j\>y^2 \, dr
\\
&<  (n-1) \int (y^{\prime})^2 \,
dr - (n-1) \int
\frac{1}{4} \Big( \frac{1
+ 4\mu ^2\chi_{[a, b]}
}{r^2}
\Big) y^2 \, dr
\\
&
= (n-1) \int (y^{\prime})^2  +
y^{\prime}{}^{\prime} y \, dr
\\
&= (n-1) (y^{\prime}y) \Big|
_{r_0}^{r_1} = 0,
\end{split}
\end{equation*}
where $I_j =
I(X_j, X_j)$,  all integrands are
evaluated at
$r$, all sums range from $j = 1$
to $j = n-1$, and all integrals are
taken from
$r_0$ to $r_1$.  This verifies
(\ref{e:sum}).

Negativity of a sum implies
negativity of at least one term,
so  (\ref{e:sum}) implies that for
some $j_0$,
$I(X_{j_0}, X_{j_0}) < 0$, and
so  by Jacobi's Theorem there
exists a point
$\gamma (r)$ with
$r_0 < r <  r_1$ which is
conjugate to
$\gamma (r_0)$.
\end{proof}

\begin{Rmk} It is
straightforward to estimate the
diameter of $M$ as follows.  All
points of $M$ lie inside the
geodesic ball at
$O$ of radius $r_1$, such that the
Sturm-Liouville solution
$y(r)$ above has its second zero
at
$r_1$.
   Here is how this reads in two
cases.

\underline{Case 1.}
$r_1 \in  (a, b]$,   a trivial
situation.  The root $r_1$ occurs
when
$
\tan (\mu  \ln (r/a)) = -\mu \ln a
$. If $a = 1$ then this gives
$$ D \leq 2 e^{\pi /\mu }.
$$

\underline{Case 2.}   $y(r) > 0$
for $a
\leq  r \leq  b$.  The second zero
of
$y(r)$  occurs at the first root of
$$
\alpha  + \beta \ln (r/b) = 0
$$ beyond $b$.  (As in the
theorem, we have
$\beta  < 0$.)  This is equivalent
to $r = be^{F }$ where
$$ F  = \frac{\alpha }{-\beta } =
\frac{\ln a
\cos \theta  +  ( \sin \theta )/\mu
}{\mu
\ln a \sin \theta  - \cos \theta }
\qquad \theta = \mu \ln (b/a),
$$ and thus,
$$ D \leq  2be^F.
$$ When $a =1$ we get
$$ D \leq  2be^{-\tan(\mu \ln
b)/\mu }.
$$ If the curvature hypotheses are
valid at all origins $O$ then the
factors $2$ in these diameter
estimates are superfluous.
   Note that as $\mu $ decreases to
$\lambda  $, the second diameter
estimate  tends to
$+\infty$.
\end{Rmk}

\section{Logarithmic Kick}
\label{s:log}
Let $e_k < r_0 \leq a < b$ be
given, where $e_k$ is the
$k^{\textrm{th}}$ superpower
of
$e$, $\ln ^k(e_k) = 0$.  Define
$\lambda  =
\lambda _k(r_0, a, b)$ as the
smallest positive solution of the
equation
$$
\cot (\lambda (\ln ^kb - \ln
^ka)) = \lambda (\ln ^ka -
\ln^kr_0).
$$
  Define
$$
F_k(r, \mu ) =
\frac{1}{4}
\Big( \frac{1}{r^2} +
\frac{1}{(r \ln r)^2}+
\dots +
\frac{1 + 4\mu ^2}{(r \ln (r)
\cdots
\ln ^k(r))^2} \Big).
$$

\begin{LogThm}
\label{t:logthm} Assume that for
all $r  \geq  r_0$ we have
$
\operatorname{Ric}(r) \geq
(n-1)F_k(r,0)
$
and that for all $r \in  [a, b]$ we
have
$$
\operatorname{Ric}(r) >
(n-1)F_k(r,\lambda ),
$$
  where $\lambda  = \lambda
_k(r_0, a, b)$ as above. Then
$M$ is compact.
\end{LogThm}

When $\mu  > 0$, the general
solution to the Boju-Funar
equation
\begin{equation}
\begin{split}
\label{e:BF}
  y^{\prime}{}^{\prime} + F_k(r,
\mu ) y = 0
\end{split}
\end{equation}
  is of the form
$$
\Phi _k(r)
\big( A \cos (\mu \ln ^kr) + B
\sin (\mu \ln ^kr) \big),
$$
where $A$, $B$ are constants
and
$$
\Phi _k(r) = \big( r \ln (r) \cdots
\ln ^{k-1}(r)\big)^{1/2}.
$$
When $\mu  = 0$ the solution
degenerates to
$$
\Phi _k(r)\big( A + B \ln ^k(r)
\big).
$$
  See \cite{BF}.

\begin{Lemm}
\label{l:BF}
If $y(r)$ solves the Boju-Funar
equation (\ref{e:BF}) with initial
conditions
$y(r_0) = 0$ and
$y^{\prime}(r_0) > 0$ and if
$\mu  > \lambda _k(r_0, a, b)$
then $y(r) = 0$ for some $r >
r_0$.
\end{Lemm}

\begin{proof}
Linear rescaling of the
$r$-variable is invalid in the
logarithmic context, but still the
proof is similar to that of
Lemma~\ref{l:SL}.  The
solution is
$$
y(r) \, = \, \Phi _k(r)
\begin{cases}
    \ln^k(r) -
\ln^k(r_0) &\textrm{ if } r_0
\leq  r
\leq  a
  \\
    A \cos (\mu \ln ^kr) + B
\sin (\mu \ln ^kr) &\textrm{ if }  a
\leq r
\leq b \\
     \alpha +
\beta  \ln ^kr  &\textrm{ if } b
\leq r <
\infty .
\end{cases}
$$
Matching values   at $a$, gives
\begin{equation}
\begin{split}
\label{e:ata}
  \ln ^k(a) - \ln^k(r_0) =
A \cos (\mu \ln ^ka) + B
\sin (\mu \ln ^ka)
\end{split}
\end{equation}
after canceling the common
factor $\Phi _k(r)$.  Matching
derivative values gives
\begin{equation}
\begin{split}
\label{e:derata}
  &\Phi
_k^{\prime}(a)
[\ln^k(a) - \ln^k(r_0)]
+
\Phi _k(a)(\ln^k)^{\prime}(a)
\\
&=
\Phi _k^{\prime}(a)\big(
A \cos (\mu \ln ^ka) + B
\sin (\mu \ln ^ka)\big)
\\
&+
\Phi _k(a) \big( -A\sin(\mu
\ln^ka) +
B\cos (\mu \ln^ka)\big)\mu \,
(\ln^k)^{\prime}(a).
\end{split}
\end{equation}
Plugging in (\ref{e:ata}),
discarding the equal terms, and
then canceling the common
factor
$\Phi _k(a)(\ln^k)^{\prime}(a)$
gives
$$
\frac{1}{\mu } =
   -A\sin(\mu
\ln^ka) +
B\cos (\mu \ln^ka) .
$$
 From this and (\ref{e:ata}) it
follows that
\begin{equation}
\begin{split}
\label{e:AB}
  A &= \cos (\mu \ln ^ka)[\ln^ka
- \ln^kr_0]- \frac{1}{\mu }
\sin (\mu  \ln^ka)
\\
B &=
\sin(\mu \ln^ka)[\ln^ka -
\ln^kr_0] + \frac{1}{\mu }
\cos(\mu \ln^ka).
\end{split}
\end{equation}
Similarly at $r = b$ we have
$$
A \cos (\mu \ln ^kb) + B
\sin (\mu \ln ^kb)
=
\alpha + \beta \ln^kb,
$$
and through more canceling we
get
$$
-A\sin(\mu \ln^kb) +
B\cos(\mu  \ln^kb) =
\frac{\beta }{\mu }.
$$
Combined with (\ref{e:AB}),
this gives
\begin{equation*}
\begin{split}
\beta  &=
\mu ( -  c(a)s(b) +
s(a)c(b)) [\ln^ka
- \ln^kr_0] + s(a)s(b) + c(a)c(b)
\end{split}
\end{equation*}
where $c(a) = \cos (\mu
\ln^ka)$,
$s(a) = \sin(\mu \ln^ka)$, etc.
But then
$$
\beta = \mu \sin (\mu (\ln^ka -
\ln^kb))[\ln^ka - \ln^kr_0]
+
\cos(\mu (\ln^ka + \ln^kb)).
$$
As in the proof of
Lemma~\ref{l:SL}, it is fair to
assume that $\mu  - \lambda $ is
small.  This ensures that $\beta
< 0$, and therefore that $y(r)$
vanishes at some $r > r_0$.
\end{proof}

\begin{proof}[\bf Proof of the
Logarithmic Kick Theorem]
Using Lemma~\ref{l:BF} in
place of Lemma~\ref{l:SL}, the
proof is the same as in the linear
case.
\end{proof}

\section{Comparison of Results}

    How does our kick  assumption
compare to that in
\cite{CGT}?  It is
different and slightly weaker.
Take, for instance $a = e$ and
$b = e^{\ell +1}$.  The
corresponding
$\lambda  = \lambda (\ell)$
tends to $0$ as
$\ell
\rightarrow
\infty$.   The   curvature
assumption in \cite{CGT} is
$$
\operatorname{Ric}(r) >
\frac{n-1}{4}\Big( \frac{1 +
4\nu ^2 }{r^2}   \Big)
$$
     for some $\nu  > 0$ and all
$r
\geq  1$.  If
$\ell$ is large then  $\lambda
(\ell)
<
\nu
$,      and for all
$r \in  [e, e^{\ell+1}]$ we have
$$
\operatorname{Ric}(r) >
\frac{n-1}{4}\Big( \frac{1+
4\lambda (\ell)^2}{r^2}  \Big),
$$ which is the hypothesis of
the Linear Kick Theorem with $a
= e$, $b = e^{\ell+1}$.  Similar
remarks are valid in the
logarithmic context.

\section{The Boundary Theorem}
\label{s:boundary}
    In this section we prove our
second main result: manifolds
with more curvature than a
Sturm-Liouville bifurcator are
compact, while those with less
curvature are noncompact. More
precisely, we assume that
$M$ is a complete
$n$-dimensional  Riemann
manifold with positive Ricci
curvature, and that $O  \in M$ is
a fixed origin.  We   fix an
SL-bifurcator
$b(r)$, i.e., a continuous
function
$b : [0,
\infty] \rightarrow  (0, \infty)$
such that the solution of
$$
w^{\prime}{}^{\prime} +
b(r)w = 0 \quad \quad
w(0) = 0
\quad w^{\prime}(0) = 1
$$
is monotone and bounded.
    Then we assert
	\begin{itemize}

    \item[(a)] $M$ is noncompact
if  for all
$r
\geq 0$,
$$
\sup \{ \operatorname{Ric}_p :
p = \exp_O(v) \textrm{ and }
\abs{v} = r\}
\leq b(r).
$$

    \item[(b)]  $M$ is compact if
$\operatorname{Ric}(r)$
exceeds $b(r)$.
\end{itemize}

\begin{proof}[\bf Proof of (a)]
This is trivial.  An
SL-bifurcator has
$$
\liminf_{r\rightarrow \infty} b(r)
= 0,
$$
whereas, the infimum of the
Ricci curvature on a compact
manifold of positive Ricci
curvature is positive.
\end{proof}

\begin{proof}[\bf Proof of (b)]
As in the proof of the Kick
Theorem, it suffices to show
that if $c(r)$ exceeds the
SL-bifurcator $b(r)$ then the
solution
$y(r)$ of the ODE
$$
y^{\prime}{}^{\prime} + c(r)y
= 0 \quad  \quad y(0) = 0 \quad
y^{\prime}(0) = 1
$$
has a second zero.  We use
Picone's Formula from
\cite{Ince}, page 226, to
compare $y(r)$ and the
corresponding  SL-bifurcator's
solution $w(r)$.
In our case, the   formula
reads:
\begin{equation*}
\begin{split}
&\frac{w(x)}{y(x)}
[w^{\prime}(x)y(x)
- w(x)
y^{\prime}(x)] \Big|_0^r
\\
= &
\int_0^r(c(x) - b(x))w(x)^2\,
dx +
\int_0^r\Big(\frac{w^{\prime}(x)
y(x)
- w(x)y^{\prime}(x)}{y(x)}
\Big)^2
\, dx,
\end{split}
\end{equation*}
provided that $y(x) > 0$ on the
interval $0 < x \leq  r$.   At $x =
0$, the contribution
$[w^{\prime}y -
wy^{\prime}]$ drops out since
$w(x)/y(x) \rightarrow 1$ as
$x \rightarrow 0$.  Thus,
$$
\frac{w(r)}{y(r)}[w^{\prime}(r)y(r)
- w(r)y^{\prime}(r)] = I(r)
$$
where $I(r)$ is the sum of the
two integrals.  Rearranging the
formula gives
\begin{equation}
\begin{split}
\label{e:y0y}
  \frac{y^{\prime}(r)}{y(r)}
=
\frac{w^{\prime}(r)}{w(r)}
-\frac{I(r)}{w^2(r)}.
\end{split}
\end{equation}
The r.h.s. of (\ref{e:y0y})
converges to a negative number
or to $-\infty$ as $r \rightarrow
\infty$.  For
$$
\lim_{r\rightarrow \infty} w(r)
>  0 \textrm{ and }
\lim_{r\rightarrow
\infty}w^{\prime}(r) = 0.
$$
Therefore there exists an  $r$
such that
$$
\frac{w^{\prime}(r)}{w(r)}
- I(r)  < 0.
$$
It follows that $y(x)$ has a
second zero.  For if $y(x)$
remains positive on $(0, r]$ then
$y ^{\prime}(r) < 0$, and
concavity implies that the graph
$y = y(x)$  subsequently crosses
the
$x$-axis.
\end{proof}

\begin{Rmk}
In \cite{Abresch} it is shown
that an SL-bifurcator satisfies
\begin{itemize}

   \item[(a)]
$\displaystyle \int_0^{\infty} r
b(r)\, dr < \infty$

   \item[(b)] All solutions of
$y^{\prime}{}^{\prime} + b(r)y
= 0$ have limit derivatives as $r
\rightarrow \infty$.

   \item[(c)]  Each solution of
$y^{\prime"}{}^{\prime} +
b(r)y = 0$ independent from
$w(r)$ diverges to $\pm
\infty$ as $r \rightarrow
\infty$.
\end{itemize}
\end{Rmk}

    \section{Planar Curves}
\label{s:planar}
    It should be easy to  descend
from higher dimensions to the
simple one dimensional case.
Unfortunately the Ricci
condition does not make much
sense. To overcome this we use
the extrinsic curvature of
curves.   The corollary to the
following result is a classification
of smooth, complete, planar
curves with non-vanishing
curvature.

We say that a smooth function
$h : [0,
\infty) \rightarrow
\rR^2$ is a
\textbf{hook} if
$\abs{h^{\prime}(s)} = 1$ for all
$s
\in  [0, \infty)$, $h$ in an
embedding (i.e., it is a
homeomorphism from $[0,
\infty)$ to the $h$-image
equipped with the inherited
topology), and the curvature
$\kappa(s)$ vanishes nowhere.
Examples of   hooks are a half
parabola and a half spiral of
infinite arclength.

\begin{HookThm}
\label{t:curves}
Under ambient
homeomorphisms of
$\rR^2$, every hook is
equivalent to one of the following
seven curves, the first four of
which are bounded:
\begin{enumerate}

\item An inward spiral to a point.
(This is ambiently
homeomorphic to a half closed
segment.)

\item

An inward spiral to a segment.

\item

An inward spiral to a circle.

\item

An outward spiral to a circle.

\item

An outward spiral to infinity.
(This is ambiently
homeomorphic to a ray.)

\item

An outward spiral to a line.  (A
line is a circle through infinity in
the
$2$-sphere.)

\item

An outward spiral to the union of
two parallel lines.
\end{enumerate}
\end{HookThm}

\begin{Cor}
\label{c:hook}
Under ambient
homeomorphisms, every
complete curve $C$ smoothly
embedded in
$\rR^2$ with
non-vanishing curvature is
equivalent to one the following
fourteen curves: the circle, the
straight line,  the twelve
combinations of the bounded and
unbounded hooks.
\end{Cor}

The proof of the Hook
Theorem  and its
corollary are left to the reader.

\begin{Rmk} For curves other
than the line, no change in the
previous corollary occurs if we
further require that the
curvature of
$C$ exceed a given positive
continuous function
$\kappa_0(s)$ with
$\dlim_{\abs{s}\rightarrow
\infty}\kappa_0(s) = 0$.  In
particular, the curvature of $C$
can be made to exceed the
curvature of a parabola.
\end{Rmk}

Despite the preceding remark,
we want to say that the
parabola lies at the boundary
between embedded and
non-embedded curves.  We
express this in two ways.

\begin{Thm}
\label{t:parabdry}
There is an arc $\kappa_t$  in the
space of curvature functions
  such that if $P_t$ is a planar
curve with curvature
  $\kappa_t$ then
\begin{itemize}

   \item[(a)]
If $t _2> t_1$ then
$\kappa_{t_2}$
exceeds $\kappa_{t_1}$.

   \item[(b)] If $t < 0$ then
$P_t$   is embedded.

   \item[(c)]  If $t > 0$ then
$P_t$ is immersed, not
embedded.

   \item[(d)]  $P_0$ is a parabola.

\item[(e)]
If $t \mapsto
\widetilde{\kappa}_t$ is an arc
of curvature functions that
approximates $t \mapsto
\kappa_t$ and $\widetilde{P}_t$
has curvature
$\widetilde{\kappa}_t$ then the
arc $t \mapsto \widetilde{P}_t$
continues to pass from
embedded curves  to
non-embedded immersed
curves at a single
time-parameter near $t = 0$.

\end{itemize}
\end{Thm}

\begin{Thm}
\label{t:paramost}
A parabola is the most curved
planar curve that is complete
and embedded as a closed
subset of the plane.
\end{Thm}

\begin{proof}[\bf Proof of
Theorem~\ref{t:parabdry}]
This is trivial.   Fix a parabola
$P_0$.  Its curvature function is
$\kappa_0$.  Fix a smooth bump
function $\beta  :
\rR \rightarrow
[0,1]$ with support $[-1,1]$, and
set
$$
\kappa_t(s) = \kappa_0(s) +
t\beta (s).
$$
It is easy to check that this arc
satisfies our assertions.
\end{proof}

To prove
Theorem~\ref{t:paramost} we
use the following lemma from
sophomore calculus.

\begin{Lemm}
\label{l:DG}
Is $C$ is a smooth curve in the
plane then the angle its tangent
turns is the integral of its
curvature.
\end{Lemm}

\begin{proof}
Parameterize $C$ by arclength,
$C(s) = (x(s), y(s))$.  The
angle turned by the tangent is
$$
\theta (s) = \angle
(C^{\prime}(0), C^{\prime}(s))
= \arctan
\Big(\frac{y^{\prime}(s)}
{x^{\prime}(s)}\Big).
$$
Its derivative is
$$
\theta ^{\prime}(s) =
\frac{(y{^{\prime}}^{\prime}x^{\prime}
- y^{\prime}
x^{\prime}{}^{\prime})/(x^{\prime})^
2}{1 +
(y^{\prime}/x^{\prime})^2}
=
(x^{\prime}{}^{\prime},
y^{\prime}{}^{\prime})
{\boldmath{\cdot}}(-y^{\prime},
x^{\prime}),
$$
since $\abs{C^{\prime}} = 1$.
For the same reason,
$C^{\prime}{}^{\prime}$ is
perpendicular  to $C^{\prime}$,
and so
$\abs{\theta ^{\prime}} =
\abs{(x^{\prime}{}^{\prime},
y^{\prime}{}^{\prime})} =
\kappa$, and the result
follows.
\end{proof}

\begin{figure}
\centering \psfig{figure = 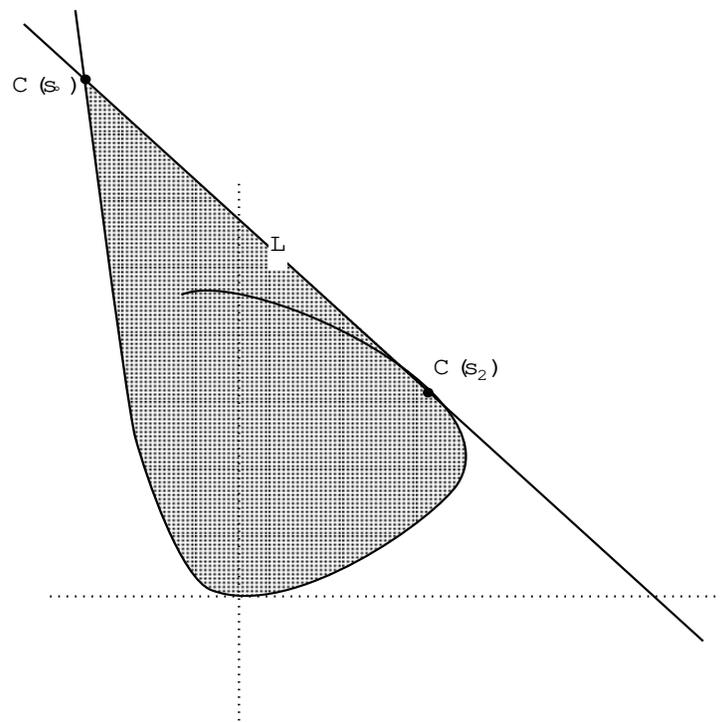,height
=95mm,width=95mm} \caption{\small {\rm  $C(s_2, \infty)$ is
trapped. }}\label{f:Case2}
\end{figure}

\begin{proof}[\bf Proof of
Theorem~\ref{t:paramost}]  The
strategy resembles that of the
Poincar\'e-Bendixson
Theorem.   Let
$P$ be a parabola
$y = kx^2$ with $k > 0$.  We
parameterize $P$ by arclength
with $P(0) = (0, 0)$.  Let $C$
be a more curved curve.  Its
curvature $\kappa(s)$ exceeds
the curvature $\kappa_P(s)$ of
the parabola, $s$ being
arclength.

By Corollary~\ref{c:hook}, $C$
is ambiently homeomorphic to a
straight line.  It has no spirals.
After translation and rotation,
$C(0) = (0, 0)$ and
$C^{\prime}(0) = (1, 0)$.
Since $\kappa$ exceeds
$\kappa_P$, there is an $s_0$
such that
$$
\kappa(s_0) > \kappa_P(s_0).
$$
We can assume $s_0> 0$.
By Lemma~\ref{l:DG},
$
\int_{0}^{\infty} \kappa_P(s) \,
ds = \pi /2,
$
so
$$
\int_0^{\infty} \kappa(s)\, ds >
\pi /2.
$$

It follows that for some first
$s_1 > 0$ we have  the turning
angle
$
\theta (s_1) = \pi /2
$.  This means that $x(s_1)$ is a
local maximum of $x(s)$ and
$C^{\prime}(s_1)$ points
vertically upward.  For
some slightly greater
$s_2$, $C^{\prime}(s_2)$
points up and leftward.  The line
$L$ through $C(s_2)$ in the
direction $C^{\prime}(s_2)$
therefore meets the positive
$y$-axis.   Since
$C$ has no spirals, $C[s_2,
\infty)$  cannot be confined to
the compact part of the plane
bounded by the
$y$-axis, $L$, and $C[0, s_2]$.
Since $C$ does not cross itself,
   it must cross the positive
$y$-axis, say at $C(s_+) = (0,
y_+)$.  See
Figure~\ref{f:crossy}.

\begin{figure}
\centering \psfig{figure =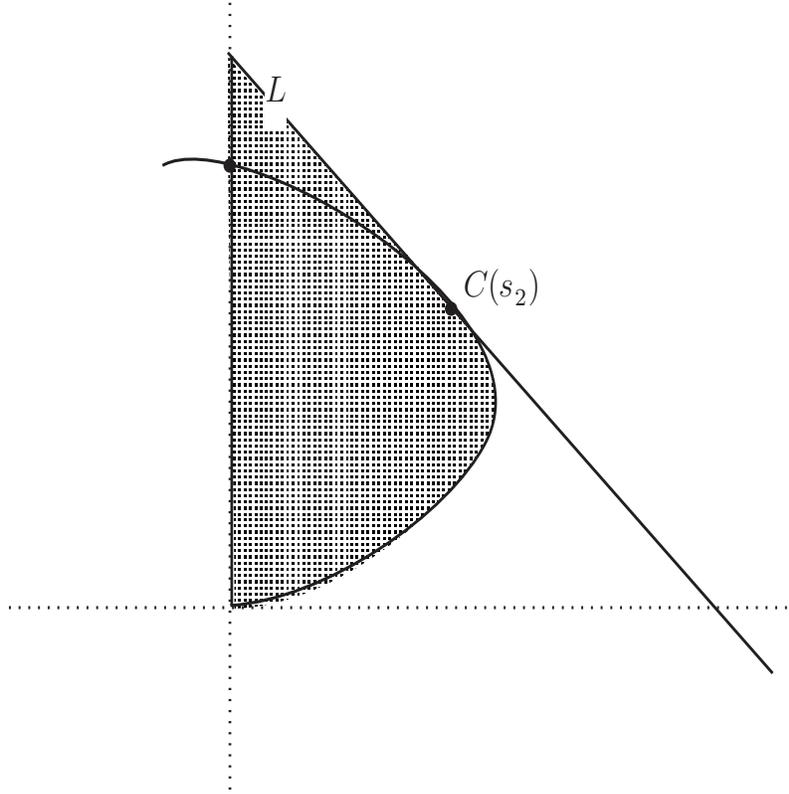,height =105mm,width=105mm} 
\caption{\small {\rm $L$ and
$C[s_2, \infty)$ cross the
$y$-axis.}} \label{f:crossy}
\end{figure}


\underline{Case  1.} For some
$s \leq  0$, $\kappa(s) >
\kappa_P(s)$.  Then $C(-\infty,
0)$ also crosses the positive
$y$-axis, say at  $C(s_-) =
(0, y_- )$ with  $0 < y_+ <
y_-$.  Then
$C(s_+, \infty) $ is
trapped inside the curve $C[s_-,
s_+]
\cup 0 \times [y_+, y_-]$, giving
a spiral, contrary to the
hypothesis on
$C$. See Figure~\ref{f:Case1}.
(If $y _-< y_+$ the roles of
$C(s_+, \infty)$ and $C(-\infty,
s_-)$ are reversed.)

\begin{figure}
\centering \psfig{figure =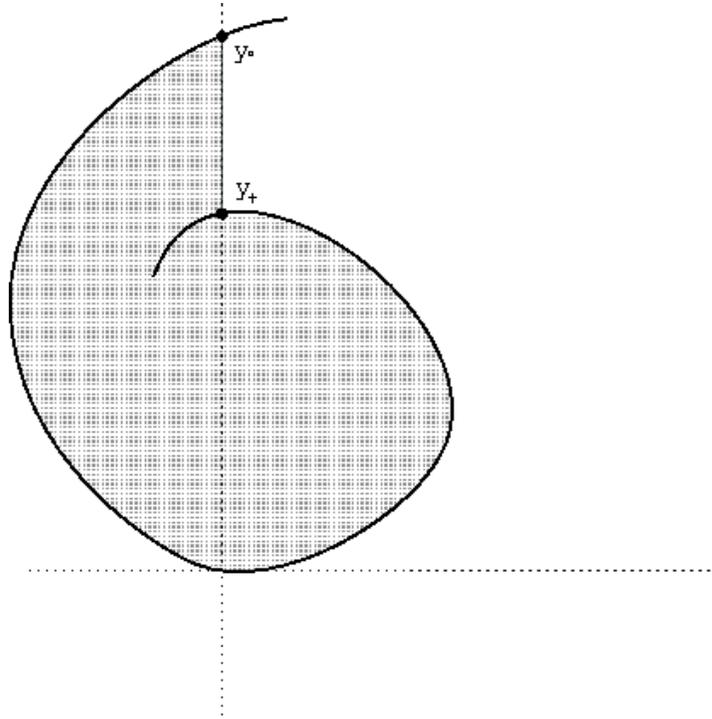,height
=95mm,width=95mm} \caption{\small {\rm  $C(s_+, \infty) $
is trapped.}}\label{f:Case1}
\end{figure}

\underline{Case  2.}   For all $s
\leq 0$, we have $\kappa(s) =
\kappa_P(s)$.  Then
$C(-\infty, 0)$ is the parabolic
arc
$y = kx^2$ with
$x < 0$.  The line $L$
crosses it, say at $C(s_-)$,
and the curve
$C(s_2, \infty)$  is confined by
  $L \cup C[s_-,s_2]$, contrary
to the fact that $C$ is ambiently
homeomorphic to a line.  See
Figure~\ref{f:Case2}.
\end{proof}

\section{Embedded Surfaces}
\label{s:embedded}

In this section we analyze
surfaces   embedded in
${\rR}^3$.
We prove two things:

\begin{Thm}
\label{t:poidmost}
Among surfaces of revolution
which are  complete
and embedded, the paraboloid is
most curved.
\end{Thm}

\begin{Thm}
\label{t:capped}
A surface of revolution that
approximates the capped
cylinder has curvature which is
an SL-bifurcator.
\end{Thm}
\begin{proof}[ \bf Proof of
Theorem~\ref{t:poidmost}]
Let $M$ be a surface of
revolution whose curvature
exceeds the curvature of the
paraboloid.  The curvature of its
profile curve exceeds the
curvature of a parabola, so by
Theorem~\ref{t:paramost} it is
not embedded.
\end{proof}

\begin{Rmk}
It seems probable that there is a
better theorem along these
lines.  It would remove the
hypothesis that $M$ is a surface
of revolution.  The reason is
that there are no surfaces
embedded in ${\rR}^3$
whose curvature decays like
that of a logarithmic (or worse)
spiral.
\end{Rmk}

\begin{proof}[\bf Proof of
Theorem~\ref{t:capped}]
Consider the surface of
revolution $M$ generated by the
  profile curve
$$
z = \frac{1}{1-\rho}
$$
in cylindrical $(\rho , \theta ,
z)$-coordinates.  The curvature
of the profile curve is $\kappa
= z^{\prime}{}^{\prime}/(1 +
(z^{\prime})^2))^{3/2}$, which is
$$
\kappa = \frac{
2}{((1-\rho )^2 +  (1-\rho
)^{-2})^{3/2}}.
$$
The curvature of $M$ at the
point $p = (\rho , \theta , z(\rho
))$ is the product of the profile
curvature and the circular
curvature, namely
$$
K(p) = \frac{  2/\rho }
{(z^{-2} +   z^2)^{3/2}}.
$$
As before, $r$ denotes the
geodesic distance from the
origin to $p$, and for large $r$,
we have $r \approx z$.  Thus
$$
K(p) \approx \frac{2}{(0 +
r^2)^{3/2}} = \frac{2}{r^3}
$$
as $r \rightarrow  \infty$.

The geodesics that emanate
from the origin are profile
curves, and since they are
asymptotically parallel as $r
\rightarrow
\infty$, we see that the solution
to the Jacobi equation
$$
w^{\prime}{}^{\prime} + b(r)w
= 0 \quad \quad w(0) = 0 \quad
w^{\prime}(0) = 1
$$
is monotone and bounded,
where $b(r)$ is the curvature at
a point
$p$ whose distance from  the
origin is $r$.  In other words,
the curvature is an
SL-bifurcator, and by the
Boundary Theorem, any surface
with more curvature is compact.
\end{proof}

\begin{Rmk}
The surface $M$ resembles a
capped cylinder.  It is
asymptotic to the cylinder $\rho
= 1$  as $r \rightarrow \infty$.
\end{Rmk}

\begin{Cor}
\label{c: }
Arbitrarily
small perturbations of the
curvature of $M$  exist
that correspond to  the topological
sphere and others that
correspond to a paraboloid-like
surface.
\end{Cor}

\begin{proof}[\bf Proof]
The curvature is an
SL-bifurcator.
\end{proof}

\begin{Rmk}
The preceding surface $M$ is
not unique.  Corresponding to any
SL-bifurcator $b(r)$ there  is a
surface of revolution whose
curvature is $b(r)$.  Its
properties  are the same as
those of $M$.
\end{Rmk}

\section{Discussion} \label{s:capped}
We have investigated  in this article
the minimium amount of energy or
curvature necessary for a complete
manifold with positive Ricci
curvature to collapse into a compact
one. The minimum amount of energy
depends on the manifold being
embedded or not. In any case, the set
of such manifolds, for which the
minimum amount of extra curvature
imposes compactness is considered
to be   the boundary between
compact and noncompact.

The following are some questions
and ideas that extend our
investigation.
\begin{enumerate}

\item  {\bf Decay versus kick.}
Above, we
   concentrated on the effect of a
curvature kick.  This amounts to
putting the Whitney topology on the
space of curvature functions, and
seeking the corresponding boundary
between compact and noncompact
manifolds.  On the other hand, one
could seek such a boundary in terms
of decay rates.  The decay rate of
curvature for the paraboloid is
$r^{-2}$ as $r \rightarrow \infty$,
while that of the smoothly capped
cylinder is $r^{-3}$.  With more
work, a capped cylinder's decay
rate can be made to be
$r^{-(2+\epsilon )}$ for any given
$\epsilon  > 0$.  (The geodesics
from the origin still are
asymptotically parallel.)  Thus, in
terms of decay rates, the paraboloid
is at the boundary of the
SL-bifurcators, although this is not
the case in terms of Whitney (kick)
perturbations.

    \item  {\bf Geometry at
infinity.} Does the decay of the Ricci
curvature impose  geometry of the
manifold at infinity?
(\emph{Topologically},
$M$ is diffeomorphic to
${\rR}^n, n\leq 3$ since its Ricci
curvature is positive.  See
\cite{GromollMeyer}and  \cite{SY}.)  For example,
if the decay of curvature for a
complete, noncompact surface $M$
is asymptotically
$r^{-3}$, are its geodesics
emanating from the origin like those
of the capped cylinder?  That is, are
they  asymptotically parallel in
the sense that they stay a bounded
distance apart?  This contrasts with
the geometry of the paraboloid,
where the curvature is
asymptotically   $r^{-2}$
and the geodesics diverge.  So in
general, we ask: is the asymptotic
geometry of
$M$ governed by the asymptotic
decay of curvature?

\item

{\bf Classification of the geometry at
infinity.}  For noncompact,
complete surfaces of positive
curvature there are two extreme
possibilities for the geometry at
infinity.  The geodesics from a fixed
origin can be asymptotically
parallel, as for the capped cylinder,
or they  can diverge, as for the
paraboloid.  Naturally, the behavior
can also be of mixed type with some
sectors of geodesics becoming
asymptotically parallel, and others
diverging.  What happens under kick
perturbation in the mixed case?  See
also Question~\ref{q:conjecture}
below.

In higher dimensions the situation
becomes more complicated.  For
example, the three dimensional pure
capped cylinder has asymptotic
geometry at infinity $\rR
\times  S^2$ and is   at the boundary
between compact and noncompact
manifolds, while the pure three
dimensional  paraboloid is not at this
boundary.    What other behavior is
there?

Passing to dimension four, we could
consider products in which we take
the  product Riemann
structure.  (It is useful to remember
that since we are dealing with Ricci
curvature, the product of manifolds
of positive Ricci curvature also has
positive Ricci curvature.)  For
example,  let
$M$ be the product of a capped
cylinder and a paraboloid.
Kick perturbations do not produce
compact manifolds, but   they can
change the geometry of $M$ at
infinity from
$({\rR} \times  S^1) \times
{\rR}^2$  to $S^2 \times
{\rR}^2$.  What is the general
geometry at infinity and how does it
change under kick (Whitney)
perturbations?

   \item
{\bf Tensors.}
On a complete noncompact
manifold $M$ with Riemann
structure $g$,  the positive Ricci
curvature condition
$\operatorname{Ric} \geq
(1/4r^2)g $ can be
re-written as
$$
\operatorname{Ric}  =
\frac{1}{4r^2}g+A
$$
where $A$ is a non-negative
symmetric tensor.
In dimension $n\geq
3$,  the second Bianchi
identity  \cite{Aubin}
implies that the tensor $A$
is not   proportional to $g$.  For
if $\operatorname{Ric} = fg$
then    $f$ is constant. Positivity
of the Ricci curvature  implies
that the constant is positive,  but
Myers' Theorem then implies
that $M$ is compact. What does
the space of tensors
$A$ that decay faster than
$r^{-2}$ look like?

   \item
{\bf Convergence.}
In what sense does a sequence of
noncompact complete manifolds
converge  to a manifold that lies on
the boundary of compact manifolds?

   \item
{\bf Singularities.}
Consider a $C^1$ $g$
such that
\begin{equation*}
\begin{split}
g&=dr^2+\alpha^2 r\sin^2(\nu
\ln(r) + \phi) d\theta^2  \textrm{ for }
a < r < b
\\
g&=dr^2+r
d\theta^2
\quad \quad \qquad \qquad \quad
\quad \textrm{ for } b
\leq  r
\end{split}
\end{equation*}
where $\phi$ and  $\alpha$ are
constants.  How can the metric be
continued smoothly in $0 < r \leq a$
so the curvature stays positive and the
singularity at $r = 0$ is minimal?
Can the singularities be classified?

   \item
{\bf SL-bifurcators in higher
dimensions.} There  are other
possibilities to define manifolds at
the boundary, using the Jacobi fields.
The  generalization of
SL-bifurcators in dimension $n\geq
3$ deals with the   matrix ODE

\begin{equation}\label{mat}
Y^{\prime}{}^{\prime} + R(s)Y = 0
\end{equation}
along geodesics $\gamma $,  where
$R(s)$  is the Ricci tensor $R(s)=
\operatorname{Ric}
(\dot{\gamma}(s),
.)\dot{\gamma}(s)$.

Consider a complete noncompact
manifold, with positive  Ricci
curvature.
Once an origin is fixed,  if the
exponential map at the origin   has no
conjugate points then the matrix
ODE has a solution $Y(t)$   such that
$Y(0) = 0$,
$Y^{\prime}(0) = I$, and
$Y(s)$ is nonsingular for $s >0$.
When $Y(s)$ has a finite limit as $s
\rightarrow \infty$ (see
\cite{Greene}, and
  \cite{Chavel} p.250),   for any
{\bf y} orthogonal to
$\dot{\gamma}(0)$,
$Y(s)${\bf y} converges to a
finite vector as $s$ tends to
infinity.

This situation corresponds  to the
case where $R(s)$ has positive
eigenvalues. In addition, the existence
of a limit at infinity imposes some
conditions on the decay of the
positive eigenfunctions.

Manifolds with such a property
could be called SL-manifolds.  When
the matrix $R(s)$ can be diagonalized
by a  matrix $P(s)$ and $P(s)$
converges to an inverstible matrix at
infinity, the matrix equation \ref{mat}
reduces to $n-1$ SL-bifurcators. In
that case, a small increase of
curvature  leads to conjugate points
and hence implies compactness. But
in general, is it true that every small
increase of curvature( $R(s)
+\epsilon Id$) of an SL-manifold
gives conjugate points and hence
implies compactness?

   \item

\label{q:conjecture}
{\bf Conjecture.}
Finally, we offer
the following simple conjecture.
Suppose that $M$ is two
dimensional, complete, and, judged
from a fixed origin $O$,  its
curvature is greater than or equal to
an SL-bifurcator.  If there is a
geodesic from $O$ that has a
conjugate point then $M$ is
compact.
\end{enumerate}

\begin{Rmk}
Clearly, if $M$ is a surface of
revolution, the conjecture is true.
But imagine a capped cylinder, and
increase its curvature on a small
open set away from the origin.  (The
perturbation depends on $\theta $.)
Does this force compactness?
\end{Rmk}


\begin{thebibliography}{20}

\bibitem{Abresch} U. Abresch,
Lower curvature bounds,
Toponogov's theorem, and bounded
topology. Annales Scientifiques de
l'\'Ecole Normale Sup\'erieure S\'er. 4,
18 no. 4 (1985), p. 651-670.

\bibitem{Aubin}T. Aubin, Some
Nonlinear Problems in Riemannian
Geometry Series : Springer
Monographs in Mathematics
1998, XVII.

\bibitem{BF} V.Boju and L. Funar, A
note on the Bonnet-Myers theorem,
Zeitsch. Analysis Anwendungen
(ZAA), 15(1996), p. 275-278.

\bibitem{Chavel} I. Chavel,
Riemmanian Geometry: A Modern
Introduction, Cambridge University
Press, Cambridge, 108, 1993.

\bibitem{CGT} J. Cheeger, M.
Gromov, and M. Taylor,  Finite
Propagation Speed, Kernel Estimates
for Functions of the Laplace
Operator and the Geometry of
Complete Riemannian Manifolds, J.
Diff. Geometry, Vol. 17, 1982, p.
15-53.

\bibitem{Kupka} B. Dekster and I.
Kupka,  Asymptotics of Curvature in
a Space of Positive Curvature,
J. Diff.
Geo., Vol. 15, 1980, p. 553-568.

\bibitem{Greene} L.W. Green,  A
theorem of E. Hopf, Mich Math. J.
1958, p. 31-34.

\bibitem{GromollMeyer} D.
Gromoll,  and W. Meyer,
On complete open manifolds  of
positive curvature.  Ann. of Math.
(2) \textbf{90} 1969, p. 75-90.



\bibitem{Gromoll} J. Cheeger and D.
Gromoll, On the structure of
complete manifolds of nonnegative
curvature, Ann. of Math. 92 (1972),
p. 413-443.


\bibitem{Ince} E. L. Ince, 
\emph{Ordinary Differential
Equations}.  Dover Publications, New
York, 1944.

\bibitem{Karcher} H. Karcher,
Riemannian comparison
constructions, in S.S.Chern Studies in
Global Geometry and Anaysis,
Studies in Math. 27, Math Assoc.
Amer. 2 nd Edn. 1989, p391-413.

\bibitem{Kobayashi} S. Kobayashi,
   On conjugate and cut loci.
1967   Studies in Global Geometry and
Analysis, Math. Assoc.
Amer.,  p. 96-122.


\bibitem{Myers} S. B. Myers,
Connections between differential
geometry and topology, Duke Math.
J. 1 (1935). p. 376-391.

\bibitem{SY} R. Schoen S.T. Yau, Complete three dimensional manifolds with positive Ricci curvature and scalar curvature, 1982 Seminar on Differential Geometry  . Math. Studies, 102 (1982), 209-228

\end{thebibliography}
\end{document}